\documentclass[12pt]{article}

\usepackage{amsfonts}
\usepackage[utf8]{inputenc}

\newtheorem{theorem}{Theorem}
\newtheorem{example}{Example}
\newtheorem{lemma}{Lemma}
\newtheorem{corollary}{Corollary}

\newtheorem{proposition}{Proposition}
\newtheorem{assumption}{Assumption}
\newtheorem{algorithm}{Algorithm}
\newtheorem{definition}{Definition}

\newcommand{\proof}{\bf Proof: \rm }%\nr}

\def\squareforqed{\hbox{\rlap{$\sqcap$}$\sqcup$}}
\def\qed{\ifmmode\else\unskip\quad\fi\squareforqed}

\def\ba{\begin{array}}
\def\ea{\end{array}}
\def\beann{\begin{eqnarray*}}
\def\eeann{\end{eqnarray*}}
\def\bea{\begin{eqnarray}}
\def\eea{\end{eqnarray}}
\def\beq{\begin{equation}}
\def\eeq{\end{equation}}
\def\BT{\begin{theorem}}
\def\ET{\end{theorem}}
\def\BE{\begin{example}}
\def\EE{\end{example}}
\def\BL{\begin{lemma}}
\def\EL{\end{lemma}}
\def\BP{\begin{proposition}}
\def\EP{\end{proposition}}
\def\BC{\begin{corollary}}
\def\EC{\end{corollary}}
\def\BD{\begin{definition}}
\def\ED{\end{definition}}
\def\BA{\begin{assumption}}
\def\EA{\end{assumption}}
\def\BAL{\begin{algorithm}}
\def\EAL{\end{algorithm}}
\def\Z{\mathbb{Z}}
\def\R{\mathbb{R}}

\begin{document}

\title{Can $n^d + 1$ unit right $d$-simplices cover a right $d$-simplex with shortest side $n + \epsilon$?}
\author{
Michael J.~Todd
\thanks{
School of Operations Research and Information Engineering,
Cornell University, Ithaca, NY 14853,
USA.
E-mail {\tt mjt7@cornell.edu}.
}}
\maketitle
\abstract{
In a famous short paper, Conway and Soifer show that $n^2 + 2$ equilateral triangles with edge length 1
can cover one with side $n + \epsilon$. We provide a generalization to $d$ dimensions.
}

\section{Introduction}
We denote by $e^1, \dots, e^d$ the unit coordinate vectors in $\R^d$, and by $e := \sum_j e^j$ the vector of ones.
A unit right $d$-simplex is defined to be the convex hull of $0$, $e^1$, $e^1+e^2$, ... , $e^1+e^2+\cdots+e^d$,
or any of its images under coordinate permutations and translations. A right $d$-simplex is a dilation of a unit
right $d$-simplex; if the dilation is by a factor $\alpha > 0$, its shortest side has length $\alpha$.

We are not able to answer the question in the title, but we do show that, if $\epsilon \leq \delta := (n+2)^{-1}$, then
$(n+1)^d + (n-1)^d - n^d$ suffice. (This fails for the trivial case $d = 1$; we assume implicitly
throughout that $d > 1$.) Notice that, under the transformation $x \mapsto Mx$, where
\[
M := \left[ \ba{cc} 1 & -1/2 \\ 0 & \sqrt{3}/2 \ea \right],
\]
right 2-simplices are transformed into equilateral triangles, so that our result implies that of Conway and Soifer
\cite{CS}.

We need a convenient notation for right $d$-simplices. For $v \in \R^d$ and $\pi$ a permutation of
$\{1,\dots,d\}$, we use $k(v,\pi)$ to denote the convex hull of $v$, $v+e^{\pi(1)}$, $v+e^{\pi(1)}+e^{\pi(2)}$,
\dots,$v+e$. It is easy to see that
\[
k(v,\pi) = \{x \in \R^d: 1 \geq (x-v)_{\pi(1)}\geq (x-v)_{\pi(2)}\geq \cdots \geq (x-v)_{\pi(d)} \geq 0 \}.
\]
It is well known that the set of all $k(0,\pi)$'s, as $\pi$ ranges over all permutations, 
triangulates the unit cube, while the set of all $k(v,\pi)$'s, with $v$ an integer vector and $\pi$ a permutation,
triangulates $\R^d$. See, for example, \cite{T}. These simplices are exactly the $d$-dimensional pieces when $\R^d$ is
partitioned by all hyperplanes of the form $x_j = z$ or $x_i - x_j = z$, with $z$ an integer. More relevant to our purposes,
the set of all $k(v,\pi)$'s, with $v$ an integer vector and $\pi$ a
permutation, that lie in the right $d$-simplex
\[
S^n := \{x \in \R^d: n \geq x_1 \geq x_2 \geq \cdots \geq x_d \geq 0 \},
\]
covers (indeed, triangulates) that set. In fact, $k(v,\pi)$ lies in this set iff
$v \in S^{n-1}$ and, if $v_j = v_{j+1}$, $j$ precedes $j+1$ in the permutation $\pi$. By volume considerations, there
are $n^d$ such unit right $d$-simplices.

We can also easily see
that the ``base'' of $S^n$, where $x_d$ lies between 0 and 1, can also be triangulated, by $n^d - (n-1)^d$
of these simplices, those with $v_d = 0$.

\section{The Result}

\BT
For $\delta := (n+2)^{-1}$, the right $d$-simplex
\[
S^{n+\delta} := \{ x \in \R^d: n + \delta \geq x_1 \geq x_2 \geq \cdot \geq x_d \geq 0 \},
\]
with shortest side $n + \delta$, can be covered by 
\[
(n + 1)^d + (n - 1)^d - n^d
\]
unit right $d$-simplices.
\ET

\proof
We divide $S^{n+\delta}$ into its base
\[
S^{n+\delta}_1 := \{ x \in S^{n+\delta}: 0 \leq x_d \leq 1+\delta \}
\]
and its top
\[
S^{n+\delta}_2 := \{ x \in S^{n+\delta}: x_d \geq 1+\delta \}.
\]

Note that the top can be written as
\[
S^{n+\delta}_2 = \{ x \in \R^d: n+\delta \geq x_1 \geq x_2 \geq \cdot \geq x_d \geq 1+\delta \},
\]
which is just the translation by $(1+\delta)e$ of $S^{n-1}$, and can therefore be triangulated by $(n-1)^d$
unit right $d$-simplices as above.

It remains to cover the base $S^{n+\delta}_1$ by $(n+1)^d - n^d$ unit right $d$-simplices. Note that this base
is somewhat similar to the base
\[
S'_1 := \{ x \in S^{n+1}: 0 \leq x_d \leq 1 \},
\]
which as we noted above, can be triangulated by exactly this many unit right $d$-simplices. Indeed, the base
we are interested in has its first $d-1$ components squeezed in (from $n+1$ to $n+\delta$) and its last
component stretched out (from 1 to $1+\delta$). We therefore apply an operation to the simplices in this
triangulation, roughly as the individual cloves are transformed by squeezing the head of a roasted garlic.

As we observed above, the simplices of the triangulation of $S'_1$ are those $k(v,\pi)$ where
\beq\label{eq:condvpi}
v \in S^n \cap \Z^d; \quad v_d = 0; \quad \mbox{   if } v_j = v_{j+1}, \mbox{ $j$ precedes $j+1$ in $\pi$}.
\eeq
We squeeze these simplices as follows: 
\[
(a) \mbox{ If } \pi^{-1}(d) = d, \tilde{k}(v,\pi) := k((1 - \delta)v, \pi); 
\]
\[
(b) \mbox{ if } \pi^{-1}(d) < d, \tilde{k}(v,\pi) := k((1 - \delta)v + \delta e, \pi).
\]
We need to show that every $x \in S^{n+\delta}_1$ is covered by at least one such $\tilde{k}(v,\pi)$,
where $(v,\pi)$ satisfies (\ref{eq:condvpi}).

For any such $x$, we can choose $v \in \Z^d_+$, $v_d = 0$, so that all components of 
\[
w := x - (1 - \delta)v,
\]
except possibly the last, lie between 0 and 1. We then order these components using the permutation $\pi$.
Suppose first we can choose $\pi$ so that $d$ comes last:
\beq\label{eq:orderw1}
1 \geq w_{\pi(1)} \geq \cdots \geq w_{\pi(d)} \geq 0, \quad \pi^{-1}(d) = d.
\eeq
Note that there is some choice involved for $j < d$; if $v_j > 0$ and $0 \leq w_j \leq \delta$, we can decrease
$v_j$ by 1 so that $1 - \delta \leq w_j \leq 1$ and
 then adjust $\pi$ accordingly. Then we have
\beq\label{eq:orderw2}
\mbox{ if } v_j > 0 \mbox{ for } 1 \leq j < d, w_j > \delta.
\eeq
Moreover, if there is a set of components of $w$ that are equal, we may modify $\pi$ so that
their indices appear in ascending order:
\beq\label{eq:orderw3}
\mbox{ if } w_{\pi(j)} = w_{\pi(j+1)} \mbox{ for } 1 \leq j < d, \pi(j) < \pi(j+1).
\eeq
We show that, if $v$ and $\pi$ can be chosen so that (\ref{eq:orderw1})--(\ref{eq:orderw3}) hold, then
$x$ lies in the simplex $\tilde{k}(v,\pi)$ of type (a). By the first of these conditions, it is only
necessary to check (\ref{eq:condvpi}).

First, we have $w_1 \geq 0$, so that
\[
v_1 \leq (1-\delta)^{-1} x_1 \leq (1-\delta)^{-1} (n+\delta) = n + 1.
\]
Moreover, if $v_1 = n+1$, we have equality throughout, so that $v_1 > 0$ and $w_1 = 0$, contradicting
(\ref{eq:orderw2}). Hence $v_1 \leq n$.

Next, consider the condition $v_j \geq v_{j+1}$. If $v_{j+1} = 0$, then this holds by default. If not,
then $j+1 < d$ and by (\ref{eq:orderw2}), $w_{j+1} > \delta$, so that $w_j < w_{j+1} + 1 - \delta$ and thus
\[
v_j > v_{j+1} - 1 + (1 - \delta)^{-1} (x_j - x_{j+1}) \geq v_{j+1} - 1,
\]
and we obtain $v_j \geq v_{j+1}$. 

Finally, if $v_j = v_{j+1}$, then since $x_j \geq x_{j+1}$ we have $w_j \geq w_{j+1}$. Thus $j$ precedes
$j+1$ in $\pi$, either by (\ref{eq:orderw1}) if these components are unequal, or by
(\ref{eq:orderw3}) if they are equal. This completes the verification of (\ref{eq:condvpi}), and so $x$ is
covered.

Note that, if $x_d \leq \delta$, then we can find $v$ and $\pi$ so that (\ref{eq:orderw1}) holds. Indeed, we order
the components of $w$ as above, and ensure that if $v_j > 0$ and $j < d$, then $w_j > \delta$ and $j$ precedes
$d$ in $\pi$. But if $v_j = 0$ for
$j < d$, then $x_j \geq x_d$ ensures that $w_j \geq w_d$, and thus we can arrange that $d$ comes last in
$\pi$. Thus the bottom sliver of the base is covered by simplices of type (a).

Now we assume that $x$ cannot be covered by such a simplex. Then $x_d > \delta$, and hence $x_j > \delta$
for all $j$. We can then find $v \in \Z^d_+$ with $v_d = 0$ and a permutation $\pi$ so that, with $w$ again
defined as $x - (1 - \delta) v$, we have
\beq\label{eq:orderw4}
1+\delta \geq w_{\pi(1)} \geq \cdots \geq w_{\pi(d)} \geq \delta.
\eeq
Moreover, as above, if $1 \leq w_j \leq 1+\delta$ for $j < d$, we can increase $v_j$ by 1 so that
$\delta \leq w_j \leq 2 \delta$ and then adjust $\pi$ accordingly, so that
\beq\label{eq:orderw5}
\mbox{ for } j < d, w_j < 1.
\eeq
We can also ensure that equal components of $w$ are suitably ordered, so that (\ref{eq:orderw3}) holds.

If $w_d > 1$, then because of (\ref{eq:orderw5}), $\pi^{-1}(d) = 1.$ If instead $w_d \leq 1$, then
(\ref{eq:orderw4}) and (\ref{eq:orderw5}) show
 that (\ref{eq:orderw1}) holds, so that if $\pi^{-1}(d) = d$, $x$ could be covered
by a simplex of type (a). Thus in either case, $\pi^{-1}(d) < d$, so that, if $w' := x - (1 - \delta)v - \delta e$,
\[
1 \geq w'_{\pi(1)} \geq \cdots \geq w'_{\pi(d)} \geq 0, \quad \pi^{-1}(d) < d,
\]
and $x$ will be covered by a simplex of type (b) if we can verify (\ref{eq:condvpi}).

Suppose (\ref{eq:orderw3})--(\ref{eq:orderw5}) hold. Then $w_1 \geq \delta$, so
\[
v_1 \leq (1-\delta)^{-1}x_1 - (1-\delta)^{-1}\delta < (1-\delta)^{-1}(n+\delta) = n+1,
\]
and we have $v_1 \leq n$.

Next, consider the condition $v_j \geq v_{j+1}$. If $v_{j+1} = 0$, then this holds by default. If not,
then $j+1 < d$ and by (\ref{eq:orderw4}) and (\ref{eq:orderw5}), $w_{j+1} \geq \delta$
and $w_j < 1$, so that $w_j < w_{j+1} + 1 - \delta$ and thus
\[
v_j > v_{j+1} - 1 + (1 - \delta)^{-1} (x_j - x_{j+1}) \geq v_{j+1} - 1,
\]
and we obtain $v_j \geq v_{j+1}$. The proof that if $v_j = v_{j+1}$ then $j$ precedes $j+1$ in the
permutation $\pi$ is identical to that above. 

Thus $x$ is covered either by a simplex of type (a) or one of type (b), and the theorem is proved.

\qed

 %{\bf Acknowledgement}
 %The author would like to thank Gena Samorodnitsky, Leonid Faybusovich, 
 %Rob Freund, Arkadi Nemirovskii, and Yurii
 %Nesterov for very helpful conversations.

\end{document}